\documentclass[10pt]{amsart}
\usepackage{amscd}
\usepackage{amssymb}
\usepackage{latexsym}
\usepackage{graphicx}
\usepackage[parfill]{parskip} 
\usepackage{amsmath}
\usepackage{amsfonts}
\usepackage{amssymb}
\usepackage{latexsym}
\usepackage{graphicx}
\usepackage{enumerate}
\usepackage{bbm}

\usepackage[utf8]{inputenc}

\usepackage[final]{hyperref}
\hypersetup{colorlinks=true, linkcolor=blue, anchorcolor=blue, citecolor=red, filecolor=blue, menucolor=blue, pagecolor=blue, urlcolor=blue}
\newcommand{\abs}[1]{\vert #1 \vert}

\newcommand{\norm}[1]{\left\Vert #1 \right\Vert}

\newcommand{\bignorm}[1]{\bigl\Vert #1 \bigr\Vert}

\newcommand{\C}{\mathbb{C}}

\newcommand{\R}{\mathbb{R}}

\newcommand{\angles}[1]{\langle #1 \rangle}
\newcommand{\bigangles}[1]{\big\langle #1 \big\rangle}

\DeclareMathOperator{\diag}{diag}

\newtheorem{theorem}{Theorem}

\newtheorem{lemma}{Lemma}

\theoremstyle{definition}

\theoremstyle{remark}
\newtheorem{remark}{Remark}

\numberwithin{equation}{section}
\title[ Yang-Mills-Higgs in Lorenz gauge]{Finite energy local well-posedness for the Yang-Mills-Higgs equations in Lorenz gauge}

\author{Achenef Tesfahun}
\address{Department of Mathematical Sciences\\ Norwegian University of Science and Technology\\ N-7491 Trondheim\\ Norway}
\email{tesfahun@math.ntnu.no}

\subjclass[2000]{35Q40; 35L70}

\begin{document}
\maketitle

\begin{abstract}
The Yang-Mills-Higgs equations (Y-M-H), when written relative to Lorenz gauge, become a system of nonlinear wave equations. The key bilinear terms in the resulting system turn out to be null forms--this is in light of a recent discovery of null structure by Selberg and the present author for Yang-Mills equations in Lorenz gauge. Using the null structure found and bilinear space-time estimates, we prove local well-posedness of Y-M-H in Lorenz gauge for finite energy data.
\end{abstract}

\section{Introduction}

The goal of this paper is to prove local well-posedness of Y-M-H in Lorenz gauge for finite energy data.
 This result is in a sense a generalization of a recent paper by Selberg and the present author \cite{st2013} on the Yang-Mills equations in Lorenz gauge for finite energy data.  
 
 Let $\mathcal G$ be a compact Lie group and $\mathfrak g$ its Lie algebra. For simplicity, we shall assume 
$\mathcal G=SO(n, \R)$ (the group of orthogonal matrices of determinant one) or $\mathcal G=SU(n, \C)$ (the group of unitary matrices of determinant one). Then $\mathfrak g=so(n, \R)$ (the algebra of skew symmetric matrices) or $\mathfrak g=su(n, \C)$ (the algebra of trace-free skew hermetian matrices)
with inner product 
\begin{align*}
\langle X,Y\rangle >=-\text{Tr}(X\cdot Y^\text{T}) \qquad \text{or} \qquad \langle X,Y\rangle >=-\text{Tr}(X\cdot Y^*),
\end{align*}
where Tr, T and * denote the trace, transpose and conjugate transpose of a matrix, respectively.  The matrix commutator is given by 
$$
[X, Y]=X\cdot Y-Y\cdot X.
$$

Given a $\mathfrak g$-valued 1-form $A$ on the Minkowski space-time $\R^{1+3}$, we denote by  $F=F^{(A)}$ the associated curvature $F=dA + [A,A]$. That is, given
 $$A_\alpha \colon \R^{1+3} \to \mathfrak g, $$
 we define $F_{\alpha\beta} =  F^{(A)}_{\alpha\beta}$ by
\begin{equation}\label{Curvature}
  F_{\alpha\beta} = \partial_\alpha A_\beta - \partial_\beta A_\alpha + [A_\alpha,A_\beta],
\end{equation}
where $\alpha,\beta \in \{0,1,2,3\}$.  

Then for a scalar field (also known as the \textit{Higgs field})
  $$\phi\colon \R^{1+3} \to \mathfrak g, $$
   Y-M-H consists of the Euler-Lagrange equations associated with the Lagrangian 
$$
\mathcal L =\langle F_{\alpha\beta}, F^{\alpha\beta}\rangle-\frac12\langle D_\alpha\phi, D^\alpha \phi\rangle-V(\phi),
$$
  where $V(\phi)$ is a \textit{Higgs potential} and $D_\alpha=D^{(A)}_\alpha$ denotes the \textit{covariant derivative} operator associated to $A$ given by $$D_\alpha = \partial_\alpha + [A_\alpha, \cdot].$$
  A typical potential is $$V(\phi)= \frac1{p+1}\abs{\phi}^{p+1}, \quad 2 \le p\le 5. $$ 
  
  In terms of $F=F^{(A)}$ and $\phi$,  Y-M-H then reads
  \begin{align}
 \label{YMH-F-1}
 D^\alpha{F_{\alpha\beta}}&=[ D_\beta \phi, \phi],\\
 \label{YMH-H-1}
 D^\alpha D_\alpha \phi &=\abs{\phi}^{p-1} \phi.
 \end{align} 
 Here and throughout this paper we follow the convention that repeated upper/lower indices are implicitly summed over their range. Indices are raised and lowered using the Minkowski metric $\diag(-1,1,1,1)$ on $\R^{1+3}$. Roman indices $i, j,k,\dots$ run over $1,2,3$ and Greek indices $\alpha, \beta, \gamma $ over $0,1,2,3$.  Points on  $\R^{1+3}$ are written $(x^0,x^1,x^2,x^3)$ with $t=x^0$, and $\partial^\alpha$ denotes the partial derivative with respect to $x^\alpha$ . We write $\partial_t=\partial_0$, $\nabla = (\partial_1,\partial_2,\partial_3)$, and $\partial=(\partial_t, \nabla)$.

 The total energy, at time $t$,  is given by
$$
  \mathcal E(t) = \int_{\R^3} \abs{F(t,x)}^2 +\abs{D\phi(t,x)}^2 + \frac1{p+1}\abs{\phi}^{p+1}  \, dx,
$$
where $D=(D_0, D_1, D_2, D_3)$. The energy is conserved for a smooth solution decaying sufficiently fast at spatial infinity,  i.e., $$ \mathcal E(t)= \mathcal E(0).$$ 
 
Given a sufficiently smooth function $\mathcal U: \R^{1+3}\rightarrow \mathcal G$, define
 the \textit {gauge transformations}
 \begin{equation}\label{GaugeTransform}
 \left\{
 \begin{aligned}
       A_\alpha\to A_\alpha'& =\mathcal U A_\alpha \mathcal U^{-1}-  (\partial_\alpha \mathcal U ) \mathcal U^{-1},
  \\
  \phi \to \phi' &=\mathcal U\phi \mathcal U^{-1}.
  \end{aligned}
         \right.
  \end{equation}
A simple calculation shows, denoting  $F'=F^{(A')}$ and $ D'_\alpha=D^{(A')}_\alpha$, that 
 \begin{align*}
   F' &=\mathcal U F \mathcal U^{-1}, \quad   D'_\alpha F'= \mathcal U [D_\alpha F] \mathcal U^{-1}, \quad
   D'_\alpha \phi'= \mathcal U [D_\alpha \phi] \mathcal U^{-1}.
\end{align*}  
These imply 
$$
 D'^\alpha{F'_{\alpha\beta}}=[ D'_\beta \phi', \phi'],\qquad
 D'^\alpha D'_\alpha \phi' =\abs{\phi'}^{p-1} \phi'.
$$
Thus,   \eqref{YMH-F-1}--\eqref{YMH-H-1} are invariant under the  \textit {gauge transformations} \eqref{GaugeTransform}, i.e., if $(F, \phi)$ satisfies \eqref{YMH-F-1}--\eqref{YMH-H-1}, so does  $(F', \phi')$.  A solution is therefore a representative of its equivalent class, and hence we may impose an additional gauge condition (on $A$). The most popular gauges are
\begin{enumerate}[(i)]
\item Temporal gauge:  $A_0=0$,
\item Coulomb gauge: $\partial^i A_i=0$,
\item Lorenz gauge:  $\partial^\alpha A_\alpha=0$.
\end{enumerate}

 In both temporal and Lorenz gauges, Y-M-H can be written as  a system of nonlinear wave equations whereas in Coulomb gauge it is expressed as a system of nonlinear wave equations coupled with an elliptic equation.  Eardley and Moncrief \cite{Eardley1982, Eardley1982b} proved local and global well-posedness of Y-M-H in the temporal gauge for initial data (for $A$ and $\phi$) in the Sobolev space \footnote{Here $H^s = (I-\Delta)^{-s/2} L^2(\R^3)$.} $H^s \times H^{s-1}$ with $s\ge 2$. 
  To prove well-posedness for data in $H^1 \times L^{2}$ (data in energy class), however, requires the bilinear terms to be null forms. In Coulomb gauge, following the work of Klainerman-Machedon \cite{Klainerman1995} on the Yang-Mills equations,  Keel \cite{Keel1997} showed that these terms are null forms and used this to prove global well-posedness of Y-M-H in Coloumb gauge for finite energy data.  For studies on the regularity theory of the related Yang-Mills equations refer to \cite{Klainerman1995, Klainerman1999, Krieger2009, Oh2012, Oh2012b, Segal1979, Sterbenz2007, st2013, Tao2003}.

Recently, Selberg and the present author \cite{st2013} discoverd null structure in the Yang-Mills equations in Lorenz gauge and subsequently proved local well-posedness for finite energy data. In this paper, we follow the same argument as in \cite{st2013} to show that Y-M-H in Lorenz gauge also contains null structure, and we combine this with space-time bilinear estimates to prove local well-posedness for finite energy data.  Lorenz-gauge null structure was first discovered in \cite{Selberg2010a} for the Maxwell-Dirac equations, and then for the Maxwell-Klein-Gordon equations \cite{Selberg2010b} (see also \cite{Selberg2013}).

\section{Y-M-H as a system of nonlinear wave equations}
 Expanding \eqref{YMH-F-1}--\eqref{YMH-H-1} yields 
 \begin{align*}
 \square{A_\beta}&=\partial_\beta\partial^\alpha A_\alpha-[ \partial^\alpha A_\alpha, A_\beta]-[A^\alpha, \partial_\alpha A_\beta]
 \\
 &\quad -[A^\alpha, F_{\alpha \beta}] -[\phi, \partial_\beta\phi] -[A^\alpha, [A_\alpha, A_\beta]]-[\phi, [A_\beta, \phi]],
 \\
 \square{\phi}&=-[\partial^\alpha A_\alpha,\phi]-2[A_\alpha, \partial^\alpha \phi]-[A^\alpha, [A_\alpha, \phi]]+\abs{\phi}^{p-1} \phi.
 \end{align*}
 Thus, in Lorenz gauge, these equations simplify to the wave equations
\begin{equation}
 \label{YMH-AP-1}
 \left\{
\begin{aligned} 
 \square{A_\beta}&=\Lambda_\beta(A, \partial A, F, \phi, \partial \phi),
 \\
 \square{\phi}&=\Phi(A,  \phi, \partial \phi),
  \end{aligned}
        \right.         
         \end{equation}
   where 
 \begin{equation}
 \label{LP}
 \left\{
\begin{aligned}
 \Lambda_\beta&=-[A^\alpha, \partial_\alpha A_\beta] -[A^\alpha, F_{\alpha \beta}] -[\phi, \partial_\beta\phi] 
 \\
&\quad -[A^\alpha, [A_\alpha, A_\beta]]-[\phi, [A_\beta, \phi]]
\\
\Phi&=-2[A_\alpha, \partial^\alpha \phi]-[A^\alpha, [A_\alpha, \phi]]+\abs{\phi}^{p-1} \phi.
\end{aligned}
        \right.         
         \end{equation}

  In addition, $F$ satisfies the wave equation 
    \begin{equation}\label{YMH-F-2}
 \begin{split}
   \square F_{\beta\gamma}&=-[A^\alpha,\partial_\alpha F_{\beta\gamma}] - \partial^\alpha[A_\alpha,F_{\beta\gamma}] - \left[A^\alpha,[A_\alpha,F_{\beta\gamma}]\right]
   \\
   & \quad - 2[F^{\alpha}_{{\;\;\,}\beta},F_{\gamma\alpha}]
   + 2[D_\beta \phi, D_\gamma\phi]+ [\phi, [F_{\beta\gamma}, \phi]].
   \end{split}
 \end{equation}
 Indeed,  this will follow if we apply apply $D^\alpha$ to the 
 Bianchi identity
$$
  D_\alpha F_{\beta\gamma} + D_\beta F_{\gamma\alpha} + D_\gamma F_{\alpha\beta} = 0
$$
and simplify the resulting expression using the commutation identity
$$
  D_\alpha D_\beta X - D_\beta D_\alpha X = [F_{\alpha\beta},X]
$$
and  \eqref{YMH-F-1}  (see  e.g. \cite{st2013}).

Expanding the second, fourth and fifth terms in \eqref{YMH-F-2} and also imposing the Lorenz gauge, yields 
\begin{equation}\label{YMH-F-3}
     \square F_{\beta\gamma}  = \Gamma_{\beta\gamma}(A, \partial A, F, \partial F,  \phi, \partial \phi),
\end{equation}
where
\begin{equation}\label{Gamma}
\left\{
\begin{aligned}
      \Gamma_{\beta\gamma} &= - 2[A^\alpha,\partial_\alpha F_{\beta\gamma}]
      + 2[\partial_\gamma A^\alpha, \partial_\alpha A_\beta]
      - 2[\partial_\beta A^\alpha, \partial_\alpha A_\gamma]
      \\
      &\quad + 2[\partial^\alpha A_\beta , \partial_\alpha A_\gamma]
      + 2[\partial_\beta A^\alpha, \partial_\gamma A_\alpha]+2[\partial_\beta \phi , \partial_\gamma \phi]
      \\
&\quad  - [A^\alpha,[A_\alpha,F_{\beta\gamma}]] + 2[F_{\alpha\beta},[A^\alpha,A_\gamma]] - 2[F_{\alpha\gamma},[A^\alpha,A_\beta]]\\
                        & \quad +  2[\partial_\beta \phi, [A_\gamma, \phi]]- 2[\partial_\gamma \phi, [A_\beta, \phi]] + [\phi, [F_{\beta\gamma},\phi]]
                        \\ 
    &\quad - 2[[A^\alpha,A_\beta],[A_\alpha,A_\gamma]]+2[[A_\beta, \phi],[A_\gamma, \phi]].
 \end{aligned}
        \right.
\end{equation}

Thus, in Lorenz gauge,  $(A, F, \phi)$ satisfies the following system of nonlinear wave equations:
\begin{equation}\label{Wave-AFP}
\left\{
\begin{aligned}
\square{A}&=\Lambda(A, \partial A, F, \phi, \partial \phi),
\\
      \square F&= \Gamma(A, \partial A, F, \partial F,  \phi, \partial \phi),
      \\
 \square{\phi}&=\Phi(A,  \phi, \partial \phi),
 \end{aligned}
        \right.
\end{equation}
where $$ \Lambda:=\{\Lambda_\beta\}, \quad \Gamma:=\{\Gamma_{\beta\gamma}\} \quad \text{and} \quad  \Phi$$
are given as in \eqref{LP} and \eqref{Gamma}. To ease notation, we often skip the arguments in $\Lambda$, $\Gamma$ and $\Phi$.

\section{The Cauchy problem and statement of the main result}
To pose the Cauchy problem for \eqref{Wave-AFP},  it suffices to consider initial data for $(A,\phi)$ at $t=0$:
\begin{equation}\label{Data-AP}
\left\{
\begin{aligned}
    A(0) &= a\in H^1, \qquad \partial_t A(0) = \dot a \in L^2,
        \\
    \phi(0) &=\phi_0 \in H^1, \qquad \partial_t \phi(0) = \phi_1\in L^2.
 \end{aligned}
        \right.
\end{equation}
 The initial data for $F$ and their regularities can be determined from \eqref{Data-AP} as follows:  Let   
 $$F(0)=f, \qquad \partial_t F(0) = \dot f . $$ Then using \eqref{Curvature} and \eqref{YMH-F-1}, we obtain
\begin{equation}\label{f}
\left\{
\begin{aligned}
  f_{ij} &= \partial_i a_j - \partial_j a_i + [a_i,a_j],
  \\
  f_{0i} &= \dot a_i - \partial_i a_0 + [a_0,a_i],
\\
  \dot f_{ij} &= \partial_i \dot a_j - \partial_j \dot a_i + [\dot a_i,a_j]+[ a_i, \dot a_j],
  \\
 \dot f_{0i} &= \partial^j f_{ji} +[a^\alpha, f_{\alpha i}] + [\partial_i\phi_0, \phi_0]+ [[a_i,\phi_0], \phi_0],
\end{aligned}
\right.
\end{equation}
where the first three expressions come from \eqref{Curvature} whereas the last one comes from \eqref{YMH-F-1} with $\beta=i$.

Using \eqref{f} and the regularities in \eqref{Data-AP}, one can show 
  $ f\in L^2$ and $\dot f \in  H^{-1}.$ Indeed, these follows as a consequence of the estimates 
\begin{equation}\label{SobHoldEst}
\left\{
\begin{aligned}
\norm{ uv}_{L^2}&\le \norm{u}_{L^4} \norm{v}_{L^4} \le  \norm{u}_{H^{\frac34}} \norm{v}_{H^{\frac34}}, 
\\
 \norm{ uv}_{H^{-1}}&\le \norm{uv}_{L^{\frac65}} \le \norm{u}_{L^2} \norm{v}_{L^3}\le \norm{u}_{L^2} \norm{v}_{H^{\frac12}},
 \\
 \norm{ uvw}_{H^{-1}}&\le \norm{uvw}_{L^{\frac65}} \le \norm{u}_{L^3} \norm{v}_{L^3} \norm{w}_{L^6} \le \norm{u}_{H^{\frac12}} \norm{v}_{H^{\frac12} }\norm{w}_{H^{1}},
\end{aligned}
        \right.
\end{equation}
where we used H\"older and Sobolev inequalities in 3d. In particular, using \eqref{SobHoldEst} in \eqref{f}, one can get the following bounds:
\begin{equation}\label{f-bound}
\left\{
\begin{aligned}
\norm{ f}_{L^2}&\le C(1+\norm{a}_{H^1})^2 \norm{\dot a}_{L^2} ,
\\
 \norm{ \dot{f}}_{H^{-1}}&\le C(1+\norm{a}_{H^1})^3 ( \norm{\dot a}_{L^2} + \norm{\phi_0}^2_{H^1}).
\end{aligned}
        \right.
\end{equation}

Thus, we have shown
\begin{equation}\label{Data-F}
  F(0) = f\in L^2,  \qquad \partial_t F(0) =\dot f \in H^{-1}.
\end{equation}
Moreover,  using \eqref{Data-AP}, H\"older inequality and Sobolev embedding, we can show 
$$ D\phi(0) \in L^2, \quad \phi_0\in L^{p+1} \quad \text{for}  \ 2 \le p\le 5 .$$
So for initial data \eqref{Data-AP}, we have $\mathcal E(0)< \infty$, and for this reason we say that the initial data have \textit{finite energy}. 

Note that \eqref{YMH-F-1} with $\beta=0$ imposes the (gauge invariant) constraint
\begin{equation}\label{DataConstraint}
  \partial^i f_{i0} + [a^i,f_{i0}] = [\phi_1, \phi_0]+ [[a_0,\phi_0], \phi_0].
\end{equation}
In addition, the Lorenz gauge imposes the constraint
\begin{equation}\label{DataConstraint-LG}
  \dot a_0=\partial^i a_i.
\end{equation}

In the rest of the paper, we shall study the system of nonlinear wave equations \eqref{Wave-AFP} in Lorenz gauge with initial data \eqref{Data-AP}, \eqref{Data-F} satisfying the constraints \eqref{DataConstraint}--\eqref{DataConstraint-LG}. The most difficult nonlinear terms to deal with in \eqref{Wave-AFP} are the bilinear terms in $\Lambda, \Gamma$ and $\Phi$. Let us for the moment focus on these bilinear terms, in which case, the system \eqref{Wave-AFP} is of the form 
\begin{equation}\label{Wave-AFP-G}
\left\{
\begin{aligned}
\square{A}&=\Pi(A,\partial A)+  \Pi(A,F)+  \Pi(\phi,\partial \phi),
\\
      \square F&= \Pi(\partial A,\partial A)+ \Pi(A,\partial F)+  \Pi(\partial\phi, \partial \phi),
      \\
 \square{\phi}&=\Pi(A, \partial \phi),
 \end{aligned}
        \right.
\end{equation}
where $\Pi(\dots)$ denotes a bilinear form in the given arguments.

It is well-known that a generic equation of the form 
$$
\square u= u\partial u \qquad (\text{or } \quad \square u= \partial u\partial u )
$$
can be shown to be  locally well-posed for $H^s\times H^{s-1}$  data for all $s>1$ (resp. $s>0$) by using Strichartz estimates 
(see \cite{p1993}).
Moreover, in view of the counter examples of Lindblad \cite{lbd1996} these results are sharp. However, if we replace $u\partial u$ by a null form of type $Q(\abs{\nabla}^{-1}u, v)$ ( resp. $\partial u\partial u$ by $Q(u, v)$), we can do better and prove local well-posedness for data in $H^1\times L^2$ (resp.  in $L^2\times H^{-1}$). Here $Q(u, v)$ is one of the following:
\begin{equation}\label{NullForm-Q}
\left\{
\begin{aligned}
Q_{0}(u,v)&=\partial_\alpha u \partial^\alpha v=-\partial_t u \partial_t v+\partial_i u \partial^j v,
\\
Q_{\alpha\beta}(u,v)&=\partial_\alpha u \partial_\beta v-\partial_\beta u \partial_\alpha v.
\end{aligned}
        \right.
\end{equation}

All the terms on the right hand side of \eqref{Wave-AFP-G} except $\Pi(A,F)$ and $ \Pi(\phi,\partial \phi)$ turn out to be null forms, as we shall show later. For this reason, we will lose regularity on the solution $A$ starting from finite energy data.  However, $A$ is only a potential representing the electromagnetic field $F$. The most interesting quantities here are $F$ and $\phi$, and we do not lose regularity for these quantities if we start with finite energy data. This is due to the presence of null structure in all the bilinear terms in the equations for $F$ and $\phi$.

We now state our main result.
\begin{theorem}\label{MainThm} Let $2 \le p<5$.  For initial data \eqref{Data-AP}, \eqref{Data-F} (which are \textit{finite energy data})  satisfying the constraints \eqref{DataConstraint}--\eqref{DataConstraint-LG},  there exists a time $T =T (\mathcal D_0)> 0$, where \footnote{Of course, $T$ depends also on $\norm{f}_{L^{2}}$ and $ \norm{\dot f}_{H^{-1}} $,  but in view of  \eqref{f-bound} these norms are bounded by some  power of $\mathcal{D}_0$.} 
$$
 \mathcal D_0= \norm{a}_{H^{1}} + \norm{\dot a}_{L^2}  +\norm{\phi_0}_{H^{1}} + \norm{\phi_1}_{L^2} ,
$$
and a solution $(A,F, \phi)$ on $(-T,T) \times \R^3$ of the system \eqref{Wave-AFP} in Lorenz gauge. The solution has the regularity \footnote{Here we use the notation $a-=a-\varepsilon$ for sufficiently small $0<\varepsilon \ll 1.$. Thus, we lose an $\varepsilon$ regularity for $A$  due to the lack of null structure in the bilinear terms of type $\Pi(A,F)$ and $ \Pi(\phi,\partial \phi)$ in the equation for $A$. }
\begin{align*}
  A &\in C([-T,T]; H^{1-}) \cap C^1([-T,T]; H^{0-}),
  \\
  F &\in C([-T,T]; L^2) \cap C^1([-T,T]; H^{-1}),
  \\
  \phi &\in C([-T,T]; H^1) \cap C^1([-T,T]; L^2),
\end{align*}
and is unique in a certain subspace of these spaces.
\end{theorem}

\begin{remark}
The semilinear wave equation
 $$
 \square \phi = \abs{\phi}^{p-1} \phi
 $$
 is energy subcritical for $p<5$ and  energy critical for $p=5$.  In the critical case, $p=5$,  local existence for data in $H^1\times L^2$ can be proved (see e.g. \cite{lbd1995}) using Strichartz estimates and a Banach fixed point argument in the space 
 $$
X_T= C([-T,T]; H^1) \cap C^1([-T,T]; L^2) \cap L^5_tL_x^{10} ([-T,T]\times \R^3).
 $$
  
 The function spaces used in the present paper, however, do not work for the critical case ($p=5$). So we consider only the subcritical case ($p<5$).
 
\end{remark}

The rest of the paper is organized as follows. In Section \ref{Sec-Null}, we reveal the null structure in the key bilinear terms in $\Lambda$,  $\Gamma$ and $\Phi$ using the Lorenz gauge. In Section \ref{Sec-Reduc}, we rewrite \eqref{Wave-AFP} as a first order system and reduce Theorem \ref{MainThm} to proving nonlinear estimates in $X^{s,b}$--spaces. In Sections \ref{Sec-NullEst} and \ref{Sec-NonNullEst}, we prove these nonlinear estimates.

  \section{Null structures in $\Lambda$,  $\Gamma$ and $\Phi$ }\label{Sec-Null}
In this Section, we follow \cite{st2013} and reveal the null structures in most of the bilinear terms in  $\Lambda$,  $\Gamma$ and $\Phi$ using the Lorenz gauge.

For $\mathfrak g$-valued $u,v$, define a commutator versions null forms by 
\begin{equation}\label{CommutatorNullforms}
\left\{
\begin{aligned}
  Q_0[u,v] &= [\partial_\alpha u, \partial^\alpha v] = Q_0(u,v) - Q_0(v,u),
  \\
  Q_{\alpha\beta}[u,v] &= [\partial_\alpha u, \partial_\beta v] - [\partial_\beta u, \partial_\alpha v] = Q_{\alpha\beta}(u,v) + Q_{\alpha\beta}(v,u).
\end{aligned}
\right.
\end{equation}
  Note the identity
\begin{equation}\label{NullformTrick}
  [\partial_\alpha u, \partial_\beta u]
  = \frac12 \left( [\partial_\alpha u, \partial_\beta u] - [\partial_\beta u, \partial_\alpha u] \right)
  = \frac12 Q_{\alpha\beta}[u,u].
\end{equation}
Define 
\begin{equation}\label{NewNull} 
  \mathfrak Q[u,v] = - \frac12 \varepsilon^{ijk}\varepsilon_{klm} Q_{ij}\left[R^l u^m, v \right]
  - Q_{0i}\left[R^i u_0, v \right],
\end{equation}
where $\varepsilon_{ijk}$ is the antisymmetric symbol with $\varepsilon_{123} = 1$ and $$R_i = \abs{\nabla}^{-1}\partial_i = (-\Delta)^{-1/2}\partial_i$$ are the Riesz transforms.

We then have the following Lemma which follows from \cite[ (2.6) in Lemma 1 and identity (2.7)]{st2013} after imposing the Lorenz gauge condition $\partial^\alpha A_\alpha=0$.
\begin{lemma}\label{NullLemma} \cite{st2013} Assume $A_\alpha,\psi \in \mathcal S$ with values in $\mathfrak g$. Then in Lorenz gauge,  $\partial^\alpha A_\alpha=0$, we have the identities
\begin{equation*}
\left\{
\begin{aligned}
 &  [A^\alpha, \partial_\alpha \psi ]  =
  \mathfrak Q\left[\abs{\nabla}^{-1} A,\psi\right],
  \\
   &[\partial_tA^\alpha, \partial_\alpha\psi]=  Q_{0i}\left[A^i , \psi \right].
\end{aligned}
\right.
\end{equation*}

\end{lemma}

Let us now look at the terms in $\Lambda$ ,  $\Gamma$ and $\Phi$ (see  \eqref{LP}, \eqref{Gamma}). In view of  Lemma \ref{NullLemma}, the first terms in $\Lambda$ and $\Phi$, and the first three terms in $\Gamma$ are all null forms. By the identity \eqref{NullformTrick}, the fourth, fifth and sixth terms in $\Gamma$ are identical to $2Q_0[A_\beta,A_\gamma]$, $Q_{\beta\gamma}[A^\alpha,A_\alpha]$ and  $Q_{\beta\gamma}[\phi, \phi]$, respectively.
Thus, the only bilinear terms which are not null forms are the second and third terms in $\Lambda$. 
    
    To this end, we split those terms which are null forms and those which are not as 
      \begin{equation}\label{Split-LGP}
    \Lambda=\Lambda^{(1)}+ \Lambda^{(2)}, \quad     \Gamma=\Gamma^{(1)}+ \Gamma^{(2)}, \quad     \Phi=\Phi^{(1)}+ \Phi^{(2)}, 
     \end{equation}  
   where $\Lambda^{(1)}, \Gamma^{(1)}$ and $ \Phi^{(1)}$ are the null form terms, and in view of the above comments these can be written as  
     \begin{equation}\label{AFP-Nonlin-1}
\left\{
\begin{aligned}
  \Lambda^{(1)}_\beta& =  -\mathfrak Q\left[\abs{\nabla}^{-1} A, A_\beta\right],
  \\
    \Gamma^{(1)}_{ij} &=- 2\mathfrak Q[\abs{\nabla}^{-1} A,F_{ij}]
  + 2\mathfrak Q[\abs{\nabla}^{-1} \partial_j A, A_i]
  - 2\mathfrak Q[\abs{\nabla}^{-1} \partial_i A, A_j]
  \\
  &\quad + 2Q_0[A_i , A_j]
  + Q_{ij}[A^\alpha,A_\alpha]+ Q_{ij}[\phi, \phi],
      \\
   \Gamma^{(1)}_{0i} &= - 2\mathfrak Q[\abs{\nabla}^{-1} A,F_{0i}]
  + 2\mathfrak Q[\abs{\nabla}^{-1} \partial_i A, A_0]
  - 2 Q_{0j}[A^j,A_i]
  \\
  &\quad + 2Q_0[A_0 , A_i]
  + Q_{0i}[A^\alpha,A_\alpha]+ Q_{0i}[\phi, \phi],
   \\
   \Phi^{(1)}&=  -2\mathfrak Q\left[\abs{\nabla}^{-1} A,  \phi\right]
 \end{aligned}
        \right.
 \end{equation}  
and 
 $\Lambda^{(2)}, \Gamma^{(2)}$ and $ \Phi^{(2)}$ are the terms without null structure:
 \begin{equation}\label{AFP-Nonlin-2}
\left\{
\begin{aligned}
  \Lambda^{(2)}_\beta& = -[A^\alpha, F_{\alpha \beta}] -[\phi, \partial_\beta\phi] 
  -[A^\alpha, [A_\alpha, A_\beta]]-[\phi, [A_\beta, \phi]],
  \\
    \Gamma^{(2)}_{\beta \gamma} &= - [A^\alpha,[A_\alpha,F_{\beta\gamma}]] + 2[F_{\alpha\beta},[A^\alpha,A_\gamma]] - 2[F_{\alpha\gamma},[A^\alpha,A_\beta]]\\
                        &\quad +  2[\partial_\beta \phi, [A_\gamma, \phi]]- 2[\partial_\gamma \phi, [A_\beta, \phi]] + [\phi, [F_{\beta\gamma},\phi]]
                        \\ 
    &\quad - 2[[A^\alpha,A_\beta],[A_\alpha,A_\gamma]]+2[[A_\beta, \phi],[A_\gamma, \phi]],
   \\
   \Phi^{(2)}&=-[A^\alpha, [A_\alpha, \phi]]+\abs{\phi}^{p-1} \phi.
 \end{aligned}
        \right.
 \end{equation}

\section{Reduction of Theorem \ref{MainThm} to nonlinear estimates} \label{Sec-Reduc}
In this Section, we rewrite \eqref{Wave-AFP} as a first order system and reduce Theorem \ref{MainThm} to proving nonlinear estimates in $X^{s,b}$--spaces. 
\subsection{Y-M-H as a first order system}
We subtract
$A$, $F$ and $\phi$ to each side of the equations in \eqref{Wave-AFP} to obtain
\begin{equation}\label{Wave-AFP-2}
\left\{
\begin{aligned}
(\square-1)A&=-A+\Lambda(A, \partial A, F, \phi, \partial \phi),
\\
     ( \square -1)F&=-F+ \Gamma(A, \partial A, F, \partial F,  \phi, \partial \phi),
      \\
 (\square-1)\phi&=-\phi+\Phi(A,  \phi, \partial \phi). 
 \end{aligned}
        \right.
\end{equation}
Thus, the wave operator $\square$ is now replaced by the Klein-Gordon operator $\square -1$. We do this to avoid singularity in the change of variables below, where we get \footnote{
 We use the notation $\angles{\cdot} =\sqrt{1+\abs{\cdot}^2}$. 
}
 $\angles{\nabla}^{-1} $ instead of the singular operator $\abs{\nabla}^{-1} $.  The change of variables are
  $$(A,\partial_t A,F,\partial_t F, \phi, \partial_t \phi) \to (A_+,A_-,F_+,F_-, \phi_+, \phi_-)$$  where
\begin{align*} 
 A_\pm & = \frac12 \left( A \pm \frac{1}{i\angles{\nabla}} \partial_t A \right),
\\
  F_\pm &= \frac12 \left( F \pm \frac{1}{i\angles{\nabla}} \partial_t F \right),
\\
  \phi_\pm &= \frac12 \left( \phi \pm \frac{1}{i\angles{\nabla}} \partial_t \phi \right).
\end{align*}
Equivalently,
\begin{equation}\label{Substitution}
\left\{
\begin{aligned}  
  (A, \  \partial_t A)
 & =
  \bigl( A_+ + A_-, \   i\angles{\nabla} (A_+ - A_-)),
  \\
   (F, \ \partial_t F)
 & =
  \bigl( F_+ + F_-,  \ i\angles{\nabla} (F_+ - F_-)),
  \\
  (\phi, \ \partial_t \phi)
 & =
  \bigl(   \phi_+ + \phi_-,  \ i\angles{\nabla} (\phi_+ - \phi_-)\bigr).
   \end{aligned}
        \right.
\end{equation}

Then the system \eqref{Wave-AFP-2} transforms to
 \begin{equation}\label{Wave-AFP-3}
\left\{
\begin{aligned}
  (i\partial_t \pm \angles{\nabla}) A_\pm &= \mp\frac{1}{2 \angles{\nabla}} \Lambda'(A_+, A_-,   F_+, F_-, \phi_+,  \phi_-),
  \\
  (i\partial_t \pm \angles{\nabla}) F_\pm &= \mp\frac{1}{2 \angles{\nabla}}\Gamma' (A_+, A_-, F_+, F_-, \phi_+,  \phi_-),
  \\
    (i\partial_t \pm \angles{\nabla}) \phi_\pm &= \mp\frac{1}{2 \angles{\nabla}} \Phi'(A_+, A_+,   \phi_+,  \phi_-),
\end{aligned}
\right.
\end{equation}
where 
\begin{equation}\label{LGP'}
\left\{
\begin{aligned}
   \Lambda'(A_+, A_-,   F_+, F_-, \phi_+,  \phi_-)&=-A+ \Lambda(A, \partial A, F, \phi, \partial \phi),
  \\
 \Gamma' (A_+, A_-, F_+, F_-, \phi_+,  \phi_-)&= -F+\Gamma(A, \partial A, F, \partial F, \phi, \partial \phi ),
  \\
   \Phi'(A_+, A_+,   \phi_+,  \phi_-)&= - \phi+ \Phi(A,  \phi, \partial \phi ).
\end{aligned}
\right.
\end{equation}
In the right-hand side of \eqref{LGP'} it is understood that we use the substitution \eqref{Substitution} on $(A, F, \phi)$ and the arguments of $\Lambda$,  $\Gamma$ and $ \Phi$. Recall also that $\Lambda$,  $\Gamma$ and $ \Phi$ are as in \eqref{Split-LGP}--\eqref{AFP-Nonlin-2}.

The initial data transforms to
\begin{equation}\label{AFP-DataSplit}
\left\{
\begin{aligned}
  A_\pm(0) &=  \frac12 \left( a \pm \frac{1}{i\angles{\nabla}} \dot a \right) \in H^{1},
  \\
  F_\pm(0) &=  \frac12 \left( f \pm \frac{1}{i\angles{\nabla}} \dot f \right) \in L^2,
  \\
   \phi_\pm(0) &= \frac12 \left( \phi_0\pm \frac{1}{i\angles{\nabla}} \phi_1 \right) \in H^1 .
\end{aligned}
\right.
\end{equation}
 
\subsection{Spaces used: $X^{s,b}$-Spaces and their properties}
We prove local well-posedness of \eqref{Wave-AFP-3}--\eqref{AFP-DataSplit} by iterating in the $X^{s,b}$-spaces adapted to the dispersive operators $i\partial_t \pm \angles{\nabla}$. These spaces are defined to be the completion of $\mathcal S(\R^{1+3})$ with respect to the norm
$$
  \norm{u}_{X^{s,b}_\pm} = \norm{\angles{\xi}^s \bigangles{-\tau \pm \angles{\xi}}^b \widetilde u(\tau,\xi)}_{L^2_{\tau,\xi}},
$$
where $\widetilde u(\tau,\xi) = \mathcal F_{t,x} u(\tau,\xi)$ is the space-time Fourier transform of $u(t,x)$.

 Let $X^{s,b}_\pm(S_T)$ denote the restriction space to a time slab $S_T = (-T,T) \times \R^3$ for $T>0$. We recall the fact that
 $$
 X^{s,b}_\pm(S_T) \hookrightarrow C([-T,T];H^s) \quad \text{for} \ b > \frac12.
 $$
 Moreover, it is well known that the linear initial value problem
 $$
  (i\partial_t \pm \angles{\nabla}) u = G  \in X^{s,b-1+\varepsilon}_\pm(S_T), \qquad u(0) = u_0\in H^s,
$$
for any $s\in \R, \ b > \frac12$,  $\ 0<\varepsilon \ll 1$,  
has a unique solution
satisfying
 \begin{equation}\label{LinearEst}
\norm{u}_{X^{s,b}_\pm(S_T)} \le C \left( \norm{u_0}_{H^s} + T^{\varepsilon} \norm{G}_{X^{s,b-1+\varepsilon}_\pm(S_T)} \right)
\end{equation}  
for $0<T<1$.

In addition to $X^{s,b}_\pm$, we shall also need the wave-Sobolev spaces $H^{s,b}$, defined to be the completion of $\mathcal S(\R^{1+3})$ with respect to the norm
$$
  \norm{u}_{H^{s,b}} = \norm{\angles{\xi}^s \bigangles{\abs{\tau}-  \angles{\xi}}^b \widetilde u(\tau,\xi)}_{L^2_{\tau,\xi}}.
$$
We shall make a frequent use of the relations
\begin{equation}\label{HXH}
\left\{
\begin{alignedat}{2}
  \norm{u}_{H^{s,b}} &\le \norm{u}_{X^{s,b}_\pm}& \quad &\text{if $b \ge 0$},
  \\
  \norm{u}_{X^{s,b}_\pm} &\le \norm{u}_{H^{s,b}}& \quad &\text{if $b \le 0$}.
\end{alignedat}
\right.
\end{equation}
In particular, \eqref{HXH} allows  us to pass from estimates in $X^{s,b}_\pm$ to corresponding estimates in $H^{s,b}$.

\subsection{Reduction to nonlinear estimates using iteration}
We shall do the iteration in the following spaces:
$$
A_\pm \in X^{s,b}_\pm(S_T), \quad F_\pm \in X^{0,b}_\pm(S_T), \quad \phi_\pm \in X^{1,b}_\pm(S_T), 
$$
where 
\begin{equation}
 s=1-\varepsilon, \qquad b=\frac12+ 2\varepsilon,
\end{equation}
for sufficiently small $0<\varepsilon \ll 1.$
Then  by \eqref{LinearEst} and a standard iteration argument, local well-posedness reduces to the  nonlinear estimates
\begin{equation}\label{ReducEst-AFP}
\left\{
\begin{aligned}
  \norm{ \Lambda'(A_+, A_-,   F_+, F_-, \phi_+,  \phi_-)}_{X^{s-1, b-1+\varepsilon}_\pm}
  &\lesssim N (1+N^4),
  \\
  \norm{ \Gamma'(A_+, A_-,   F_+, F_-, \phi_+,  \phi_-)}_{X^{-1,b-1+\varepsilon}_\pm}
  &\lesssim N (1+N^4),
  \\
    \norm{ \Phi'(A_+, A_-,    \phi_+,  \phi_-)}_{X^{0,b-1+\varepsilon}_\pm}
  &\lesssim N (1+N^4),
\end{aligned}
\right.
\end{equation} 
where
\begin{align*}
 N& = \sum_\pm\left( \norm{A_\pm}_{X^{s, b}_\pm} 
  + \norm{F_\pm}_{X^{0,b}_\pm} + \norm{\phi_\pm}_{X^{1,b}_\pm}\right).
\end{align*}

 The estimates for the linear terms $A$, $F$ and $\phi$ in $\Lambda'$, $ \Gamma'$ and $\Phi'$ (see \eqref{LGP'}) are trivial, and so we ignore them.  Then recalling \eqref{Split-LGP}--\eqref{AFP-Nonlin-2}, the estimates in \eqref{ReducEst-AFP} reduce to proving
\begin{equation}\label{AFP-ReducEst-1}
\left\{
\begin{aligned}
  \norm{\Lambda^{(1)}}_{H^{s-1,b-1+\varepsilon}}
  &\lesssim N (1+N^4),
  \\  
   \norm{\Gamma^{(1)}}_{H^{-1,b-1+\varepsilon}}
  &\lesssim N (1+N^4),
  \\  
    \norm{\Phi^{(1)}}_{H^{0,b-1+\varepsilon}}
  &\lesssim N (1+N^4),
\end{aligned}
\right.
\end{equation} 
and
\begin{equation}\label{AFP-ReducEst-2}
\left\{
\begin{aligned}
  \norm{\Lambda^{(2)}}_{H^{s-1,b-1+\varepsilon}}
  &\lesssim N (1+N^4),
  \\  
   \norm{\Gamma^{(2)}}_{H^{-1,b-1+\varepsilon}}
  &\lesssim N (1+N^4),
  \\  
    \norm{\Phi^{(2)}}_{H^{0,b-1+\varepsilon}}
  &\lesssim N (1+N^4),
\end{aligned}
\right.
\end{equation} 
where we have also used \eqref{HXH} to replace  $X^{s,b}_\pm$ type norms on the left hand sides by $H^{s,b}$ (we can do this since $b-1+\varepsilon<0$).  In order to ease notation, we have not written the arguments for $\Lambda^{(i)}$, $ \Gamma^{(i)}$ and $\Phi^{(i)}$ but it is understood that the arguments are there (as in \eqref{AFP-Nonlin-1}--\eqref{AFP-Nonlin-2}).

  \section{Proof of \eqref{AFP-ReducEst-1} }\label{Sec-NullEst}
  
  \subsection{Reduction of \eqref{AFP-ReducEst-1} to estimates for  $Q=Q_0$, $Q_{0i}$ or $Q_{ij}$}
 To simplify notation, for $u=u_++u_-$, we write
 $$
 \norm{u}_{X^{s,b}} =\norm{u_+}_{X_+^{s,b}} +\norm{u_-}_{X_-^{s,b}} .
 $$
Looking at the terms in  $\Lambda^{(1)}$, $ \Gamma^{(1)}$ and $\Phi^{(1)}$  (see \eqref{AFP-Nonlin-1}) and noting the fact that the Riesz transforms $R_i$ are bounded in the spaces involved, the estimates in \eqref{AFP-ReducEst-1} reduce to proving the corresponding estimates for the null forms $Q$ :
 \begin{align}
  \label{NullEst-A11}
  \norm{Q[\abs{\nabla}^{-1} A, A]}_{H^{s-1,b-1+\varepsilon}}
  &\lesssim \norm{A}_{X^{s,b}} \norm{A}_{X^{s,b}},
  \\
  \label{NullEst-F11}
  \norm{Q[\abs{\nabla}^{-1} A, F]}_{H^{-1, b-1+\varepsilon}}
  &\lesssim \norm{A}_{X^{s, b}} \norm{F}_{X^{0,b}},
  \\
  \label{NullEst-F12}
  \norm{Q[A, A]}_{H^{-1, b-1+\varepsilon}}
  &\lesssim \norm{A}_{X^{s,b}} \norm{A}_{X^{s, b}},
  \\
    \label{NullEst-F13}
  \norm{Q[\phi, \phi]}_{H^{-1,b-1+\varepsilon}}
  &\lesssim \norm{\phi}_{X^{1,b}} \norm{\phi}_{X^{1, b}},
  \\
  \label{NullEst-P11}
  \norm{Q[\abs{\nabla}^{-1} A,  \phi]}_{H^{0, b-1+\varepsilon}}
  &\lesssim \norm{A}_{X^{s,b}} \norm{\phi}_{X^{1,b}}.
    \end{align}

    \subsection{Reduction of \eqref{NullEst-A11}--\eqref{NullEst-P11} to estimates for ordinary null forms}
The matrix commutator null forms are linear combinations of the ordinary ones, in view of \eqref{CommutatorNullforms}. Since the matrix  structure plays no role in the estimates under consideration,  we reduce to estimates to the ordinary null forms for $\mathbb C$-valued functions $u$ and $v$ (as in \eqref{NullForm-Q}).

Now, substituting
\begin{align*}
u&=u_+ + u_-, \quad \partial_t u = i\angles{\nabla}(u_+ - u_-),
\\
   v&=v_+ + v_-, \quad \partial_t v = i\angles{\nabla}(v_+ - v_-),
\end{align*}  
one obtains
\begin{align*}
  Q_0(u,v)
  &=
  \sum_{\pm,\pm'} (\pm 1)(\pm' 1)  \bigl[ \angles{D} u_{\pm} \angles{D} v_{\pm'} - (\pm D^i)u_{\pm} (\pm' D_i)v_{\pm'} \bigr],
  \\
  Q_{0i}(u,v)
  &=
  \sum_{\pm,\pm'} (\pm 1)(\pm' 1) \bigl[ - \angles{D} u_{\pm} (\pm' D_i)v_{\pm'} + (\pm D_i)u_{\pm} \angles{D}v_{\pm'} \bigr],
  \\
  Q_{ij}(u,v)
  &=
  \sum_{\pm,\pm'} (\pm 1)(\pm' 1) \bigl[ - (\pm D_i)u_{\pm} (\pm' D_j)v_{\pm'} + (\pm D_j)u_{\pm} (\pm' D_i)v_{\pm'} \bigr],
\end{align*}
where
$$
  D = (D_1,D_2,D_3) = \frac{\nabla}{i}
$$
has Fourier symbol $\xi$. 
In terms of the Fourier symbols
\begin{align*}
  q_0(\xi,\eta) &= \angles{\xi}\angles{\eta} - \xi \cdot \eta,
  \\
  q_{0i}(\xi,\eta) &= - \angles{\xi} \eta_i + \xi_i \angles{\eta},
  \\
  q_{ij}(\xi,\eta) &= - \xi_i \eta_j + \xi_j \eta_i,
\end{align*}
we have
\begin{align*}
  Q_0(u,v) &= \sum_{\pm,\pm'} (\pm 1)(\pm' 1) B_{q_0(\pm \xi, \pm' \eta)}( u_{\pm}, v_{\pm'} ),
  \\
  Q_{0i}(u,v) &= \sum_{\pm,\pm'} (\pm 1)(\pm' 1) B_{q_{0i}(\pm \xi, \pm' \eta)}( u_{\pm}, v_{\pm'} ),
  \\
  Q_{ij}(u,v) &= \sum_{\pm,\pm'} (\pm 1)(\pm' 1) B_{q_{ij}(\pm \xi, \pm' \eta)}( u_{\pm},  v_{\pm'}),
\end{align*}
where for a given symbol $\sigma(\xi,\eta)$ we denote by $B_{\sigma(\xi,\eta)}(\cdot,\cdot)$ the operator defined by
$$
  \mathcal F_{t,x} \left\{ B_{\sigma(\xi,\eta)}(u,v) \right\}(\tau,\xi) = \int \sigma(\xi-\eta,\eta) \widetilde u(\tau-\lambda,\xi-\eta) \widetilde v(\lambda,\eta) \, d\lambda \, d\eta.
$$

The symbols appearing above satisfy the following estimates.
\begin{lemma}\label{Lmm-nullsyb} \cite{st2013} For all nonzero $\xi,\eta \in \R^3$,
\begin{align*}
  \abs{q_0(\xi,\eta)} &\lesssim \abs{\xi}\abs{\eta}\theta(\xi,\eta)^2 + \frac{1}{\min(\angles{\xi},\angles{\eta})},
  \\
  \abs{q_{0j}(\xi,\eta)} &\lesssim \abs{\xi}\abs{\eta}\theta(\xi,\eta) + \frac{\abs{\xi}}{\angles{\eta}} + \frac{\abs{\eta}}{\angles{\xi}},
  \\
  \abs{q_{ij}(\xi,\eta)} &\le \abs{\xi}\abs{\eta} \theta(\xi,\eta),
\end{align*}
\end{lemma}
where $\theta(\xi,\eta) = \arccos\left(\frac{\xi\cdot\eta}{\abs{\xi}\abs{\eta}}\right) \in [0,\pi]$ is the angle between $\xi$ and $\eta$.  It is this angle which quantifies the null structure in the bilinear terms. This is due to the following angle estimate, which allows us to trade in hyperbolic regularity and gain a corresponding amount of elliptic regularity.
\begin{lemma}\label{AngleEstimate} Let $\alpha,\beta,\gamma \in [0,1/2]$. Then for all pairs of signs $(\pm,\pm')$, all $\tau,\lambda \in \R$ and all nonzero $\xi,\eta \in \R^3$,
$$
  \theta(\pm\xi,\pm'\eta)
  \lesssim
  \left( \frac{\angles{\abs{\tau+\lambda}-\abs{\xi+\eta}}}{\min(\angles{\xi},\angles{\eta}} \right)^\alpha
  +
  \left( \frac{\angles{-\tau\pm\abs{\xi}}}{\min(\angles{\xi},\angles{\eta})} \right)^\beta
  +
  \left( \frac{\angles{-\lambda\pm'\abs{\eta}}}{\min(\angles{\xi},\angles{\eta})} \right)^\gamma.
$$
\end{lemma}
For a proof, see for example \cite[Lemma 2.1]{Selberg2008}.

In view of Lemma \ref{Lmm-nullsyb}, and since the norms we use only depend on the absolute value of the space-time Fourier transform, we can reduce any estimate for $Q(u,v)$ to a corresponding estimate for the three expressions
$$
  B_{\theta(\pm\xi,\pm'\eta)}(\abs{\nabla} u, \abs{\nabla} v),
  \quad
  \angles{\nabla} u \angles{\nabla}^{-1} v
  \quad \text{and} \quad
  \angles{\nabla}^{-1} u \angles{\nabla} v.
$$
Thus, \eqref{NullEst-A11}--\eqref{NullEst-P11} can be reduced to the following:
\begin{align}
  \label{NullEst1:1}
  \norm{B_{\theta(\pm\xi,\pm'\eta)}(u, v)}_{H^{s-1,b-1+\varepsilon}}
  &\lesssim \norm{u}_{X^{s,b}_\pm} \norm{v}_{X^{s-1, b}_{\pm'}},
  \\
  \label{NullEst2:1}
  \norm{B_{\theta(\pm\xi,\pm'\eta)}(u, v)}_{H^{-1, b-1+\varepsilon}}
  &\lesssim \norm{u}_{X^{s,b}_\pm} \norm{v}_{X^{-1,b}_{\pm'}},
  \\
  \label{NullEst3:1}
  \norm{B_{\theta(\pm\xi,\pm'\eta)}(u, v)}_{H^{-1, b-1+\varepsilon}}
  &\lesssim \norm{u}_{X^{s-1, b}_\pm} \norm{v}_{X^{s-1, b}_{\pm'}},
   \\
  \label{NullEst4:1}
  \norm{B_{\theta(\pm\xi,\pm'\eta)}(u, v)}_{H^{-1, b-1+\varepsilon}}
  &\lesssim \norm{u}_{X^{0,b}_\pm} \norm{v}_{X^{0,b}_{\pm'}},
  \\
  \label{NullEst5:1}
  \norm{B_{\theta(\pm\xi,\pm'\eta)}(u, v)}_{H^{0, b-1+\varepsilon}}
  &\lesssim \norm{u}_{X^{s, b}_\pm} \norm{v}_{X^{0, b}_{\pm'}},
  \\
  \label{NullEst1:2}
  \norm{uv}_{H^{s-1,0}}
  &\lesssim \norm{\abs{\nabla}u}_{H^{s-1, b}} \norm{v}_{H^{s+1,b}},
  \\
  \label{NullEst1:3}
  \norm{uv}_{H^{s-1,0}}
  &\lesssim \norm{\abs{\nabla}u}_{H^{s+1, b}} \norm{v}_{H^{s-1,b}},
  \\
  \label{NullEst1:4}
  \norm{uv}_{H^{-1,0}}
  &\lesssim \norm{\abs{\nabla}u}_{H^{s-1,b}} \norm{v}_{H^{1,b}},
  \\
  \label{NullEst1:5}
  \norm{uv}_{H^{-1,0}}
  &\lesssim \norm{\abs{\nabla}u}_{H^{s+1,b}} \norm{v}_{H^{-1,b}},
  \\
  \label{NullEst1:6}
  \norm{uv}_{H^{-1,0}}
  &\lesssim \norm{u}_{H^{s-1,b}} \norm{v}_{H^{s+1,b}},
  \\
  \label{NullEst1:7}
  \norm{uv}_{H^{-1,0}}
  &\lesssim \norm{u}_{H^{0,b}} \norm{v}_{H^{2,b}},
  \\
  \label{NullEst1:8}
  \norm{uv}_{L^2}
  &\lesssim \norm{\abs{\nabla}u}_{H^{s-1, b}} \norm{v}_{H^{2,b}},
  \\
  \label{NullEst1:9}
  \norm{uv}_{L^2}
  &\lesssim \norm{\abs{\nabla}u}_{H^{s+1, b}} \norm{v}_{H^{0,b}}.
\end{align}

\subsection{Proof of \eqref{NullEst1:1}--\eqref{NullEst1:9}}
To prove these estimates, we need the following  Theorem about product estimates in $H^{s,b}$--spaces which is due to D'Ancona, Foschi and Selberg  \cite{dfs2006}.
\begin{theorem}\label{ThmAtlas} \cite{dfs2006}. Let $s_0,s_1,s_2 \in \R$ and $b_0,b_1,b_2 \ge 0$. Assume that
\begin{equation*}
\left\{
\begin{aligned}
  \sum b_i &> \frac12,
  \\
  \sum s_i &> 2 - \sum b_i,
  \\
  \sum s_i& > \frac32 - \min_{i \neq j}(b_i+b_j),
  \\
  \sum s_i& > \frac32 - \min(b_0 + s_1 + s_2, s_0 + b_1 + s_2, s_0 + s_1 + b_2),
  \\
  \sum s_i& \ge 1,
  \\
  \min_{i \neq j}(s_i+s_j) &\ge 0,
\end{aligned}
\right.
\end{equation*} 
and that the last two inequalities are \textbf{not both equalities}. Then
$$  
\norm{uv}_{H^{-s_0,b_0}} \lesssim \norm{u}_{H^{s_1,b_1}} \norm{v}_{H^{s_2,b_2}}
$$  
holds for all $u,v \in \mathcal S(\R^{1+3})$.
  \end{theorem}
   We also need  the following null form estimates which is proved in \cite{st2013} using the angle estimate in Lemma \ref{AngleEstimate} and the product estimates in Theorem \ref{ThmAtlas} above.
 \begin{theorem}\label{NullformThm} \cite{st2013}
Let $\sigma_0,\sigma_1,\sigma_2,\beta_0,\beta_1,\beta_2 \in \R$. Assume that
\begin{equation*}
\left\{
\begin{aligned} \label{NullformThm:1}
 & 0 \le \beta_0 < \frac12 < \beta_1,\beta_2 < 1,
  \\
 & \sum \sigma_i + \beta_0 > \frac32 - (\beta_0 + \sigma_1 + \sigma_2),
  \\
 & \sum \sigma_i > \frac32 - (\sigma_0 + \beta_1 + \sigma_2),
  \\
&  \sum \sigma_i > \frac32 - (\sigma_0 + \sigma_1 + \beta_2),
  \\
  &\sum \sigma_i + \beta_0 \ge 1,
  \\
 & \min(\sigma_0 + \sigma_1, \sigma_0 + \sigma_2, \beta_0 + \sigma_1 + \sigma_2) \ge 0,
\end{aligned}
\right.
\end{equation*} 
and that the last two inequalities are \textbf{not both equalities}. Then we have the null form estimate
$$
  \norm{B_{\theta(\pm\xi,\pm'\eta)}(u,v)}_{H^{-\sigma_0,-\beta_0}}
  \lesssim
  \norm{u}_{X^{\sigma_1,\beta_1}_\pm} \norm{v}_{X^{\sigma_2,\beta_2}_{\pm'}}.
$$
\end{theorem}

Then the null form estimates \eqref{NullEst1:1}-\eqref{NullEst5:1} follow from Theorem \ref{NullformThm} and the estimates \eqref{NullEst1:6} and \eqref{NullEst1:7}  hold by Theorem \ref{ThmAtlas}. 

We remain to prove \eqref{NullEst1:2}-\eqref{NullEst1:5} and \eqref{NullEst1:8}-\eqref{NullEst1:9}.  Consider first the case where $u$ has Fourier support on $\abs{\xi} \ge 1$. Then we can replace $\abs{\nabla} u$ by $\angles{\nabla} u$. Consequently, we can replace the norms $\norm{\abs{\nabla}u}_{H^{s-1, b}} $ and $\norm{\abs{\nabla}u}_{H^{s+1, b}} $ on the rights hands of those estimates by $\norm{u}_{H^{s, b}} $ and $\norm{u}_{H^{s+2, b}} $  respectively, and the resulting new estimates will hold by Theorem \ref{ThmAtlas}.
 
 Next, consider the case where $u$ has Fourier support on $\abs{\xi} <1$. Then the frequency $\eta$ of $v$ is comparable to the output frequency $\xi+\eta$, in the sense that $\angles{\xi+\eta} \sim \angles{\eta}$. Hence we reduce to
$$
  \norm{uv}_{L^2} \lesssim \norm{\abs{\nabla}u}_{L^2} \norm{v}_{H^{0,b}}.
$$
But this follows by estimating
$$
  \norm{uv}_{L^2} \le
  \norm{u}_{L_t^2(L_x^\infty)} \norm{v}_{L_t^\infty(L_x^2)}
$$
followed by the low-frequency estimate
\begin{equation*}
  \norm{u}_{L_t^2(L_x^\infty)} \lesssim
   \bignorm{ \norm{\widehat u(t,\xi)}_{L_\xi^{1} }}_{L^2_t} \le 
   \norm{ \abs{\xi}^{-1}}_{L_\xi^2  ( \abs{\xi} <1)}\norm{  
   \norm{ \abs{\xi}  \widehat u(t,\xi) }_{ L_\xi^2}  }_{ L^2_t}\lesssim  \norm{ \abs{\nabla}u }_{L^2}
\end{equation*}
 and the estimate
\begin{equation*}
  \norm{v}_{L_t^\infty(L_x^2)} \lesssim \bignorm{ \norm{\widetilde v(\tau,\xi)}_{L_\tau^{1}} }_{L^2_\xi} \le \norm{\angles{\tau}^{-b}}_{L_\tau^{2 }} \norm{v}_{H^{0,b}}\lesssim \norm{v}_{H^{0,b}}.
\end{equation*}

\section{Proof of \eqref{AFP-ReducEst-2}}\label{Sec-NonNullEst}
In the estimate for $\Phi^{(2)}$ it suffices to consider $p=4$. 
Looking at the terms in  $\Lambda^{(2)}$, $ \Gamma^{(2)}$ and $\Phi^{(2)}$  in \eqref{AFP-Nonlin-2} and ignoring the matrix commutator structure, the estimates in \eqref{AFP-ReducEst-2} reduce to
 \begin{align}
\label{NonLinEst-A-13}
 \norm{uv}_{H^{s-1, b-1-\varepsilon}}
  &\lesssim \norm{u }_{H^{s, b}} \norm{v }_{H^{0, b}},
  \\
  \label{NonLinEst-A-14}
  \norm{ u \angles{\nabla}v}_{H^{s-1,b-1-\varepsilon}}
  &\lesssim   \norm{u }_{H^{1, b}} \norm{v }_{H^{1, b}},
  \\
  \label{NonLinEst-A-15}
  \norm{ uvw}_{H^{s-1, b-1-\varepsilon}}
  &\lesssim \norm{u }_{H^{s, b}}\norm{v }_{H^{s, b}}\norm{w }_{H^{s, b}},
  \\
    \label{NonLinEst-F-13}
  \norm{ uvw}_{H^{-1,b-1-\varepsilon}}
  &\lesssim \norm{u }_{H^{s, b}}\norm{v }_{H^{s, b}}\norm{w }_{H^{0, b}},
  \\
    \label{NonLinEst-F-14}
  \norm{ uv\angles{\nabla}w}_{H^{-1,b-1-\varepsilon}}
  &\lesssim \norm{u }_{H^{s, b}}\norm{v }_{H^{1, b}}\norm{w }_{H^{1, b}},
  \\
     \label{NonLinEst-F-15}
  \norm{ uvwz}_{H^{-1,b-1-\varepsilon}}
  &\lesssim \norm{u }_{H^{s, b}}\norm{v }_{H^{s, b}}\norm{z }_{H^{s, b}}\norm{w }_{H^{s, b}},
  \\
     \label{NonLinEst-P-13}
  \norm{ uvw}_{H^{0,b-1-\varepsilon}}
  &\lesssim \norm{u }_{H^{s, b}}\norm{v }_{H^{s, b}}\norm{w }_{H^{1, b}}.
  \end{align}
Note that estimates for some of the terms in  $\Lambda^{(2)}$, $ \Gamma^{(2)}$ and $\Phi^{(2)}$ are ignored since the estimates for those terms are either similar or weaker than the estimates we have in \eqref{NonLinEst-A-13}--\eqref{NonLinEst-P-13}.

 The estimate \eqref{NonLinEst-A-13} follows from Theorem \ref{ThmAtlas}  and \eqref{NonLinEst-A-14} reduces to  
  $$
 \norm{ u v}_{H^{s-1, b-1-\varepsilon}}
  \lesssim   \norm{u }_{H^{1, b}} \norm{v }_{H^{0, b}}  
  $$
  which is weaker than \eqref{NonLinEst-A-13}. We estimate 
  \eqref{NonLinEst-A-15} as
\begin{align*}
\norm{ uvw}_{H^{s-1, b-1-\varepsilon}}&\lesssim \norm{ uv}_{H^{\frac12, 0}}\norm{ w}_{H^{1, b}}
\\
&\lesssim \norm{ u}_{H^{s, b}}  \norm{ v}_{H^{s, b}} \norm{ w}_{H^{1, b}},
\end{align*}    
  where we used  Theorem \ref{ThmAtlas} in both steps. 
  
 We estimate \eqref{NonLinEst-F-13} using  Theorem \ref{ThmAtlas} as
\begin{align*}
\norm{ uvw}_{H^{-1, b-1-\varepsilon}}&\lesssim \norm{ uv}_{H^{\frac12, 0}}\norm{ w}_{H^{0, b}}
\\
&\lesssim \norm{ u}_{H^{s, b}}  \norm{ v}_{H^{s, b}} \norm{ w}_{H^{0, b}},
\end{align*}    
  whereas  \eqref{NonLinEst-F-14} can be reduced to estimate \eqref{NonLinEst-F-13}. We estimate \eqref{NonLinEst-F-15}  as
\begin{align*}
\norm{ uvwz}_{H^{-1, b-1-\varepsilon}}&\lesssim \norm{ uvw}_{L^2}\norm{ z}_{H^{s, b}}
\\
&\lesssim \norm{ uv}_{H^{\frac12+\varepsilon, 0}}  \norm{ w}_{H^{s, b}}
\norm{ z}_{H^{s, b}}
\\
&\lesssim \norm{ u}_{H^{s, b}}  \norm{ v}_{H^{s, b}} \norm{ w}_{H^{s, b}}
\norm{ z}_{H^{s, b}}.
\end{align*}    
    
  Finally, we estimate  \eqref{NonLinEst-P-13} as
    \begin{align*}
\norm{ uvw}_{H^{0, b-1-\varepsilon}}&\lesssim \norm{ uv}_{H^{\frac12, 0}}\norm{ w}_{H^{1, b}}
\\
&\lesssim \norm{ u}_{H^{s, b}}  \norm{ u}_{H^{s, b}} \norm{ w}_{H^{1, b}}.
\end{align*}

\subsection*{Acknowledgement}
The author would like to thank Sigmund Selberg for encouraging him to work on this problem and for his collaboration on an earlier paper concerning Yang-Mills equations.


\begin{thebibliography}{9}

\bibitem{dfs2006}
P. D'Ancona, D. Foschi, and S. Selberg, 
\emph{Atlas of products for wave-{S}obolev spaces on
  {$\mathbb{R}^{1+3}$}}, Trans. Amer. Math. Soc. \textbf{364} (2012), no.~1,
  31--63. \MR{2833576 (2012g:46049)}
  
  
\bibitem{Selberg2010a}
Piero D'Ancona, Damiano Foschi, and Sigmund Selberg, \emph{Null structure and
  almost optimal local well-posedness of the {M}axwell-{D}irac system}, Amer.
  J. Math. \textbf{132} (2010), no.~3, 771--839. \MR{2666908 (2011h:81117)}


\bibitem{Eardley1982}
D.~M. Eardley and V. Moncrief, \emph{The global existence of
  {Y}ang-{M}ills-{H}iggs fields in {$4$}-dimensional {M}inkowski space. {I}.
  {L}ocal existence and smoothness properties}, Comm. Math. Phys. \textbf{83}
  (1982), no.~2, 171--191. \MR{649158 (83e:35106a)}

\bibitem{Eardley1982b}
\bysame, \emph{The global existence of {Y}ang-{M}ills-{H}iggs fields in
  {$4$}-dimensional {M}inkowski space. {II}. {C}ompletion of proof}, Comm.
  Math. Phys. \textbf{83} (1982), no.~2, 193--212. \MR{649159 (83e:35106b)}

\bibitem{Keel1997}
M. Keel, \emph{Global existence for critical power {Y}ang-{M}ills-{H}iggs
  equations in {${\bf R}^{3+1}$}}, Comm. Partial Differential Equations
  \textbf{22} (1997), no.~7-8, 1161--1225. \MR{1466313 (99b:58227)}

\bibitem{Klainerman1995}
\bysame, \emph{Finite energy solutions of the {Y}ang-{M}ills equations in
  {$\bold R^{3+1}$}}, Ann. of Math. (2) \textbf{142} (1995), no.~1, 39--119.
  \MR{1338675 (96i:58167)}

\bibitem{Klainerman1999}
S. Klainerman and D. Tataru, \emph{On the optimal local regularity for
  {Y}ang-{M}ills equations in {${\bf R}^{4+1}$}}, J. Amer. Math. Soc.
  \textbf{12} (1999), no.~1, 93--116. \MR{1626261 (2000c:58052)}

\bibitem{Krieger2009}
J.~Krieger, W.~Schlag, and D.~Tataru, \emph{Renormalization and blow up for the
  critical {Y}ang-{M}ills problem}, Adv. Math. \textbf{221} (2009), no.~5,
  1445--1521. \MR{2522426 (2010h:58023)}

\bibitem{lbd1996}
H. Lindblad, \emph{Counterexamples to local existence for semi-linear wave
  equations}, Amer. J. Math. \textbf{118} (1996), no.~1, 1--16. \MR{1375301
  (97b:35124)}


\bibitem{lbd1995}
H. Lindblad, C. Sogge \emph{On existence and scattering with minimal regularity for semilinear wave
  equations}, Journal of Functional Analysis \textbf{130}, 357--426 (1995).
  
\bibitem{Oh2012}
S.J. Oh, \emph{Gauge choice for the {Y}ang-{M}ills equations using the
  {Y}ang-{M}ills heat flow and local well-posedness in ${H}^{1}$}, preprint
  (2012), arXiv:1210.1558, 2012.

\bibitem{Oh2012b}
\bysame, \emph{Finite energy global well-posedness of the {Y}ang-{M}ills
  equations on $\mathbb{R}^{1+3}$: An approach using the {Y}ang-{M}ills heat
  flow}, preprint (2012), arXiv:1210.1557, 2013.

\bibitem{p1993}
G. Ponce and T. C. Sideris, \emph{Local regularity of nonlinear wave
  equations in three space dimensions}, Comm. Partial Differential Equations
  \textbf{18} (1993), no.~1-2, 169--177. \MR{1211729 (95a:35092)}
  
  
\bibitem{Segal1979}
Irving Segal, \emph{The {C}auchy problem for the {Y}ang-{M}ills equations}, J.
  Funct. Anal. \textbf{33} (1979), no.~2, 175--194. \MR{546505 (84a:58025)} 
  
\bibitem{Selberg2010b}
Sigmund Selberg and Achenef Tesfahun, \emph{Finite-energy global well-posedness
  of the {M}axwell-{K}lein-{G}ordon system in {L}orenz gauge}, Comm. Partial
  Differential Equations \textbf{35} (2010), no.~6, 1029--1057. \MR{2753627
  (2011m:35371)}

\bibitem{Selberg2013}
\bysame, \emph{Global well-posedness of the {C}hern-{S}imons-{H}iggs equations
  with finite energy}, Discrete Contin. Dyn. Syst. \textbf{33} (2013), no.~6,
  2531--2546. \MR{3007698}


\bibitem{st2013}
\bysame, \emph{Null structure and local well-posedness in the energy class for the Yang-Mills equations in Lorenz gauge.},
http://arxiv.org/pdf/1309.1977.pdf

\bibitem{Selberg2008}
Sigmund Selberg, \emph{Anisotropic bilinear {$L^2$} estimates related to the
  3{D} wave equation}, Int. Math. Res. Not. IMRN (2008), Art. ID rnn 107, 63.
  \MR{2439535 (2010i:35216)}


\bibitem{Sterbenz2007}
J. Sterbenz, \emph{Global regularity and scattering for general non-linear
  wave equations. {II}. {$(4+1)$} dimensional {Y}ang-{M}ills equations in the
  {L}orentz gauge}, Amer. J. Math. \textbf{129} (2007), no.~3, 611--664.
  \MR{2325100 (2008e:58037)}


\bibitem{Tao2003}
T.  Tao, \emph{Local well-posedness of the {Y}ang-{M}ills equation in the
  temporal gauge below the energy norm}, J. Differential Equations \textbf{189}
  (2003), no.~2, 366--382. \MR{1964470 (2003m:58016)}





\end{thebibliography}
\end{document}